\documentclass[12pt]{amsart}
\usepackage{amssymb,amsmath}
\headsep=1cm \numberwithin{equation}{section}
\theoremstyle{definition}

\theoremstyle{remark}

\newcommand{\fg}{finitely generated }

\newcommand{\Hom}{\operatorname{Hom}}

\newcommand{\Spec}{\operatorname{Spec}}
\newcommand{\Min}{\operatorname{Min}}
\newcommand{\Max}{\operatorname{Max}}
\newcommand{\Ass}{\operatorname{Ass}}
\newcommand{\Assh}{\operatorname{Assh}}

\newcommand{\Supp}{\operatorname{Supp}}
\newcommand{\rank}{\operatorname{rank}}

\newcommand{\fm}{\mathfrak{m}}
\newcommand{\fp}{\mathfrak{p}}
\newcommand{\fq}{\mathfrak{q}}
\newcommand{\fa}{\mathfrak{a}}

\newcommand{\fr}{\mathfrak{r}}

\newcommand{\im}{\operatorname{Im}}

\newcommand{\depth}{\operatorname{depth}}
\newcommand{\Coker}{\operatorname{Coker}}
\newcommand{\Ker}{\operatorname{Ker}}
\newcommand{\heit}{\operatorname{ht}}
\newcommand{\Id}{\operatorname{Id}}

\begin{document}

\author[Mohammad T. Dibaei]{Mohammad T. Dibaei }
\title[A study of Cousin complexes through the dualizing complexes]
{A study of Cousin complexes through the  dualizing complexes}
 \keywords{Cousin complexes, dualizing complexes, Gorenstein modules.\\
 \indent This research is supported in part by MIM Grant P82--104.\\
 \indent Address: Mosaheb Institute of Mathematics, Teacher Training University,
  599 Taleghani Avenue, 19165 Tehran, Iran.\\
 \indent E--mail address: dibaeimt@ipm.ir .}
\address {Institute for Studies in Theoretical Physics and
Mathematics, P.O.Box 19395-5746, Tehran, Iran.}

\begin{abstract}For the Cousin complex of certain modules, we
investigate finiteness of cohomology modules, local duality
property and injectivity of its terms. The existence of canonical
modules of Noetherian non-local rings and   the Cousin complexes
of them with respect to the height filtration are discussed .
\end{abstract}

\maketitle

\section*{Introduction}
This Paper is the continuation of [DT1] and [DT2]. We have seen in
[DT2] that if $M$ is a {\fg} module over a local ring $A$ which
possesses the fundamental dualizing complex $I^{\bullet}$,
 then, under certain conditions on $M$, ${\Hom}_{A}(M,I^{\bullet}$)
represents the Cousin complex of the module
$H^{{\dim}A-{\dim}_{A}(M)}({\Hom}_{A}(M,I^{\bullet}$)),
 the $( {\dim} {A} - {{\dim}_{A}} (M))$-th cohomology module of the complex
 ${\Hom}_{A}($M$,I^{\bullet}$), with respect to an appropriate
 filtration of ${\Spec}(A)$; and that we can reconstruct the Cousin
 complex of the module $M$ by means of the fundamental dualizing
 complex (see the proof of [DT2, Lemma 3.1]). In
 section 2, we pursue our expectation that the Cousin complexes of such modules will inherit some properties of
 the dualizing complex of the ring itself. We will show that if
 $({A}, {\fm})$ is a Noetherian local ring ( not necessarily
 possessing a dualizing complex) such that all of its formal
 fibres are Cohen-Macaulay rings, $M$ is a finitely generated $A$-module
 which satisfies the condition $({S}_{2})$ of Serre and
 ${\Min}_{\widehat{A}}(\widehat{M})={\Ass}_{\widehat{A}}(\widehat{M})$,
  then all the cohomology modules
 of $C_{A}(M)$, the Cousin complex of $M$ with respect to the
 $M$-height filtration, are finitely generated $A$-modules (a result proved also,
 under different hypotheses, by T. Kawasaki in [K, Theorem 1.1]), and
 also they satisfy a local duality property which is analogous to
 that of the Grothendieck local duality. Here, $\widehat{M}$ denotes the
 completion of $M$ with respect to the ${\fm}$-adic topology. We
 present, in section 3, a number of applications which come out of these results
 and those of [DT1] and [DT2] .\\
  \indent In the remainder of the paper we study the Cousin complex
 of certain modules over Noetherian non-local ring $R$. In section 4
 we recall the notion of canonical modules for such a ring $R$ and
 prove the existence of them when $R$ possesses a dualizing
 complex and satisfies $({S}_{2})$. As a result we present a partial
  generalization of [BH, Proposition 3.3.18]. In section 5, we generalize  [DT1, Corollary
3.4] for non-local case
 and give a characterization for the Cousin
 complex of a canonical module w.r.t. the height filtration to be
 a dualizing complex.
 Finally, we give  an explicit description for all indecomposables
 injective modules which improves [DT1, Corollary 3.3 ] .

 \section{  Preliminaries}
 Throughout $A$ is a Noetherian local ring of dimension $d$ with
 the maximal ideal ${\fm}$, and $M$ is a finitely generated
 $A$--module of ${\dim}_{A}(M)=s$. A finitely generated $A$-module
 $K_{M}$ (if it exists) is called the $canonical$ $module$ of $M$ if
 $K_{M}\otimes_{A}\widehat{A}\cong{\Hom}_{A}({{H}}^{s}_{\fm}(M),
 {E}({A}/{\fm}))$,
 where ${{H}}^{s}_{\fm}(M)$ is the s--th local cohomology module
 of $M$ w.r.t. ${\fm}$,  and ${E}({A}/{\fp})$  is the injective
 envelope of the $A$--module ${A}/{\fp}$ with
 ${\fp}\in{\Spec}({A})$. The canonical module of $M$ (if exists)
 is unique up to isomorphism (see [HK, Lemma 5.8] ).
 \subsection{Some remarks}If $A$ possesses a dualizing complex,  then
 it possesses the fundamental dualizing complex
\begin{center} ${I^{\bullet} : 0 \longrightarrow I^{0}
\overset{\delta^{0}}{\longrightarrow} {I}^{1}\overset{\delta^{1}}{
\longrightarrow} \cdots \overset{\delta^{d-1}}{\longrightarrow}
{I}^{d} \longrightarrow0}$,
\end{center}
\noindent which we call `` the dualizing complex " (see [H]), with
the
 following properties :\break
   \indent(i) for each $i\geq0$,  $H^{i}(I^{\bullet})$, the $i$-th cohomology
 module of $I^{\bullet}$, is finitely generated.\break
   \indent(ii) $I^{i} =\underset{{\fp}\in{{\Spec}(A)}, {\dim}(A/{\fp})=d-i}{\bigoplus}
   E(A/{\fp})$,              $i=0, 1, ..., d$.\break
   \noindent If $A$ possesses the dualizing complex $I^{\bullet}$, then the
   module $K_{M}:= H^{d-s}({\Hom}_{A}(M, I^{\bullet}))$ is the
   canonical module of $M$. If $K_{M}$ is the canonical module
   of $M$, it is easy to see that $\widehat{(K_{M})}\cong{K_{\widehat{M}}}$
   is the canonical module of $\widehat{M}$, as
   $\widehat{A}$-module.\break
For the module $M$, we set ${\Min}_{A}(M)$ to denote the set of
all minimal elements of ${\Supp}_{A}(M)$, and
\begin{center}${\Assh}_{A}(M) =\{{\fp}\in{\Supp}_{A}(M) :
{\dim}_{A}(A/{\fp})={\dim}_{A}(M)\}$.
\end{center}
\noindent Also $M$ is said to satisfy $(S_{n})$ if
${\depth}_{{A}_{\fp}}({M}_{\fp}) \geq {\min}\{n,
{\heit}_{M}({\fp})\}$ for all ${\fp} \in {\Supp}_{A}(M)$.
\noindent A filtration of ${\Spec}(A)$ is a descending sequence
$\mathcal{F} = (F_{i})_{i\geq0}$   of subsets of ${\Spec}(A)$, so
that, \begin{center} $F_{0} \supseteq F_{1} \supseteq \cdots
\supseteq F_{i}\supseteq\cdots, $ \end{center} with the property
that, for each $i\geq0$, each member of $\partial F_{i}
=F_{i}-F_{i+1}$ is a minimal member of $F_{i}$, with respect to
inclusion . We say that $\mathcal{F}$ admits $M$ if
${\Supp}_{A}(M)\subseteq F_{0}$. Suppose that $\mathcal{F}$ is a
filtration of ${\Spec}(A)$ that admits $M$. The Cousin complex
$C(\mathcal{F}, M)$ for $M$ with respect to $\mathcal{F}$ has the
form \begin{center} $0 \overset{d^{-2}}{\longrightarrow}
M\overset{d^{-1}}{\longrightarrow} M^{0} \overset{d^{0}}
{\longrightarrow} M^{1} \longrightarrow \cdots \longrightarrow
 M^{n-1} \overset{d^{n-1}}{\longrightarrow} M^{n} \longrightarrow \cdots$
 \end{center}
 with $M^{n} = \bigoplus_{{\fp} \in {\partial{F}_{n}}}({\Coker}
 d^{n-2})_{\fp}$ for all $n\geq 0$, and with differentiation
 $d^{n}$, as recalled in [T].\\
\indent Set ${\mathcal{H}}_{M} = (H_{i})_{i\geq0}$ to be the
 $M$-height filtration of ${\Spec}(A)$,  i.e.  $H_{i} = \{{\fp} \in
 {\Supp}_{A}(M) : {\heit}_{M}({\fp}) \geq i\}$. We denote the Cousin
 complex of $M$ with respect to ${\mathcal{H}}_{M}$ by $C_{A}(M)$.

\indent Set $\bar{A} = A/0:_{A}M$. Then $M$ has a natural
structure as $\bar{A}$-module. It is straightforward to see that
each term of the complex $C_{A}(M)$ has a natural $\bar{A}$-module
structure and each  differentiation of $C_{A}(M)$ is an
$\bar{A}$--homomorphism. Moreover, it is straightforward to see
that :
\subsection{Lemma} {\it If M  is a finitely generated A-module and
} $\bar{A}:=A/0:_{A}M$, {\it then there exists an isomorphism of
complexes} $C_{A}(M)\cong C_{\bar{A}}(M)$.\\

 \indent The following lemma will be used later.\\
  \subsection{Lemma} [P, Theorem 3.5] {\it Suppose that all formal
  fibres of A are Cohen-Macaulay. If M is a finitely generated
  A-module, then there is a morphism of complexes } $u^{\bullet}:
  C_{A}(M)\bigotimes_{A}\widehat{A}\longrightarrow C_{\widehat{A}}(\widehat{M})$
  {\it which is a monomorphism. Moreover the quotient complex}
  $Q^{\bullet}$, {\it in the exact sequence }
  \begin{center} $0 \longrightarrow C_{A}(M)\bigotimes_{A}\widehat{A}
   \overset{u^{\bullet}} \longrightarrow C_{\widehat{A}}(\widehat{M})
    \longrightarrow Q^{\bullet}\longrightarrow 0$,
  \end{center}
   is an exact complex, so that, for each $i\geq0$,
   there exists an $\widehat{A}$--isomorphism $H^{i}(C_{A}(M))
   \bigotimes_{A}\widehat{A} \cong H^{i}(C_{\widehat{A}}(\widehat{M}))$.

\subsection{Convention} For a complex $C^{\bullet}: 0\longrightarrow
C^{-1}\overset{\theta^{-1}}\longrightarrow
C^{0}\overset{\theta^{0}}\longrightarrow C^{1}\overset{\theta^{1}}
\longrightarrow\cdots$ of $A$ --module and $A$--homomorphisms, we
denote
\begin{center} $C': 0 \longrightarrow
C^{0}\overset{\theta^{0}} \longrightarrow C^{1}
\overset{\theta^{1}} \longrightarrow \cdots$ and $({C'})^{\ast}:0
\longrightarrow H^{0}(C')\longrightarrow C^{0}\longrightarrow
C^{1}\longrightarrow \cdots$.\end{center}
\section{Some properties of Cousin complexes} In this section we
establish some properties of certain complexes by means of
dualizing complexes. First we show that these Cousin complexes
have finitely generated cohomologies.
\subsection{Theorem} {\it Let A be a ring with Cohen-Macaulay
formal fibres. Assume that M satisfies $(S_{2})$ and }
${\Min}_{\widehat{A}}(\widehat{M})={\Assh}_{\widehat{A}}(\widehat{M})$.
{\it Then }
${C_{A}(M)'}$ {\it has finitely generated cohomology modules}.\\

\noindent {\it Proof}. Since $M$ satisfies $(S_{2})$, the Cousin
complex $C_{A}(M)$ is exact at $M$ and $M^{0}$ (see [SSc, Example
4.4]). Thus $H^{0}(C_{A}(M)')= M$. So it is enough to prove that
$H^{i}(C_{A}(M)')$ is finitely generated for all $i>0$. Note that,
for $i>0$, we have $H^{i}(C_{A}(M)')=H^{i}(C_{A}(M))$. \indent By
1.3, we have $H^{i}(C_{A}(M))\bigotimes_{A}\widehat{A}\cong
H^{i}(C_{\widehat{A}}(\widehat{M}))$. Therefore
$C_{\widehat{A}}(\widehat{M})$ is also exact at $\widehat{M}$ and
$(\widehat{M})^{0}$ ; so that $\widehat{M}$ satisfies $(S_{2})$ as
$\widehat{A}$--module. Since ${\Min}_{\widehat{A}}(\widehat{M}) =
{\Assh}_{\widehat{A}}(\widehat{M})$, by [DT2, Theorem 3.2], all
cohomology modules $H^{i}(C_{\widehat{A}}(\widehat{M}))$ are
finitely generated $\widehat{A}$--module. Now, by [M, Exercise
7.3], the claim follows. $\Box$
\subsection{Corollary} {\it Assume
that the ring $A$ satisfies $(S_{2})$ and all formal fibres of
$A$ are Cohen-Macaulay. Then $C_{A}(A)$, the Cousin complex of
$A$, has finitely generated cohomology modules}.\\

\noindent {\it Proof.}  By [M, Theorem 23.9], $\widehat{A}$
satisfies
  $(S_{2})$ and
thus ${\Min}(\widehat{A})={\Assh}(\widehat{A})$ (see[DT1, Remark
1.3]). $\Box$\\
 \indent For a ring $A$ and a property $P$, the $P$
locus of $A$ is defined to be the set $P(A)= \{
{\fp}\in{\Spec}(A): P$ holds for $A_{\fp}\}$. We show that the
$(S_{n})$ locus of any $(S_{2})$ local ring with Cohen-Macaulay
formal fibres is an open subset of ${\Spec}(A)$ for all $n\geq 2$
.
\subsection{Corollary} {\it If A satisfies} $(S_{2})$ {\it and all formal
fibres of } $A\longrightarrow \widehat{A}$ {\it are Cohen-Macaulay
, then for each} $n\geq 0$, $S_{n}(A)$ {\it is an open subset of}
${\Spec}(A)$, {\it in the Zariski topology. In particular, CM(A)
is an open subset of} ${\Spec}(A)$.\\

\noindent {\it Proof.}  It follows that $\widehat{A}$ is
$(S_{2})$. We assume that $n\geq 3$. Set $U_{i}=
{\Spec}(A)-{\Supp}_{A}(H^{i}(C_{A}(A))), 1\leq i\leq n-2$. Each
$U_{i}$ is an open subset of ${\Spec}A$, because
${\Supp}_{A}(H^{i}(C_{A}(A)))=V(0:_{A}H^{i}(C_{A}(A)))$ by 2.2.
Set $W=\cap_{i=1}^{n-2}U_{i}$. We show that $S_{n}(A)=W$. Let
${\fp} \in S_{n}(A)$; so that $A_{\fp}$ is $(S_{n})$. Thus,  by
[SSc, Example 4.4], $H^{n}(C_{A_{\fp}}(A_{\fp})) = 0$  for $1\leq
i\leq n-2$. Therefore,  by [S1, Theorem 3.5], we have that
${\fp}\in U_{n}$ for all $i$, $1\leq i\leq n-2$; that is ${\fp}\in
W$. In a
similar way,  we have $W\subseteq S_{n}(A)$.$\square$ \\

 \indent Next, we state a local duality property for
the Cousin complexes of certain modules.
\subsection{Theorem}(Local duality for certain Cousin
complexes). {\it Assume that all formal fibres of A are
Cohen-Macaulay, M satisfies $(S_{2})$, and that
${\Min}_{\widehat{A}}(\widehat{M})={\Assh}_{\widehat{A}}(\widehat{M})$.
Then, for each $i\geq 0$, $D_{A}H^{i}(C_{A}(M)') \cong
H_{\widehat{{\fm}}}^{s-i}(K_{\widehat{M}})$, where $D_{A}
:={\Hom}_{A}(- , E(A/{\fm}))$. Moreover, if M admits a canonical
module, then the completion signs on the right hand side of the
above isomorphism can be removed}.\\

 \noindent {\it Proof}.  Set
$\bar{A}=A/0:_{A}M$ and
$\bar{\widehat{A}}=\widehat{A}/0:_{\widehat{A}} \widehat{M}$. It
is straightforward to see that $\bar{\widehat{A}}$ and
$\widehat{\bar{A}}$ are isomorphic rings. Let $J^{\bullet}$ be the
dualizing complex for $\widehat{A}$ and assume that
$I^{\bullet}={\Hom}_{\widehat{A}}(\bar{\widehat{A}}, J^{\bullet})$
such that $I^{0}={\Hom}_{\widehat{A}}(\bar{\widehat{A}},
J^{d-s})$. Hence $I^{\bullet}$ is the dualizing complex for
$\bar{\widehat{A}}$. As seen in the proof of 2.1, $\widehat{M}$
satisfies $(S_{2})$ as $\widehat{A}$--module. It is easy to see
that $\widehat{M}$ also satisfies $(S_{2})$ as
$\bar{\widehat{A}}$-module. Since
${\dim}_{\bar{\widehat{A}}}(\widehat{M})={\dim}(\bar{\widehat{A}})$,
we have, by the proof of [DT2, Lemma 3.1], the isomorphism of
complexes
$$C_{\bar{\widehat{A}}}(\widehat{M})'\cong
{\Hom}_{\bar{\widehat{A}}}(K_{\widehat{M}}, I^{\bullet}).$$

\indent Therefore, by 2.1 and [B-ZS, Corollary 2.5], we have

\begin{equation} D_{\bar{\widehat{A}}}H^{i}(C_{\bar{\widehat{A}}}(\widehat{M})') \cong
H_{\bar{\widehat{{\fm}}}}^{s-i}(K_{\widehat{M}})\tag{1}\end{equation}
for all $i\geq 0$.

On the other hand each formal fibre of $\bar{A}$ is also a formal
fibre of $A$ and $\bar{\widehat{A}}\cong\widehat{\bar{A}}$. Hence,
from 1.3, we have
\begin{equation} H^{i}(C_{\widehat{\bar{A}}}(\widehat{M})')\cong H^{i}(C_
{\bar{A}}(M)')\otimes_{\bar{A}}\widehat{\bar{A}},
\tag{2}\end{equation} for all $i>0$. From (1) and (2), we obtain
\begin{equation}
{\Hom}_{\widehat{\bar{A}}}(H^{i}(C_{\bar{A}}(M)')\otimes_{\bar{A}}
\widehat{\bar{A}},
E_{\widehat{\bar{A}}}(\widehat{\bar{A}}/\widehat{\bar{{\fm}}}))
\cong H_{\widehat{\bar{{\fm}}}}^{s-i}(K_{\widehat{M}}).\tag{3}
\end{equation}The left hand side of (3) is
isomorphic to $${\Hom}_{\bar{A}}(H^{i}(C_{\bar{A}}(M)'),
{\Hom}_{\widehat{\bar{A}}}(\widehat{\bar{A}},
E_{\widehat{\bar{A}}}(\widehat{\bar{A}}/\widehat{\bar{{\fm}}})))$$
which, in turn, is isomorphic to
${\Hom}_{\bar{A}}(H^{i}(C_{\bar{A}}(M)'),
E_{\widehat{\bar{A}}}(\widehat{\bar{A}}/\widehat{\bar{{\fm}}}))$.
Thus, we have from (3), the isomorphism
\begin{equation}{\Hom}_{\bar{A}}(H^{i}(C_{\bar{A}}(M)'),
 E_{\bar{A}}({\bar{A}}/{\bar{{\fm}}})) \cong
 H_{\bar{\widehat{{\fm}}}}^{s-i}(K_{\widehat{M}}).\tag{4}
\end{equation}

Assume that $N$ is an $A$, $\bar{A}$--bimodule such that, for $a
\in A$ and $x \in N$, $ax=\bar{a}x$, where $-:A
\longrightarrow\bar{A}$ is the natural map. Then we have

$$\begin{array}{llll}
{\Hom}_{\bar{A}}(N, E_{\bar{A}}(\bar{A}/\bar{{\fm}}))
&\cong{\Hom}_{\bar{A}}(N, {\Hom}_{A} (\bar{A}, E(A/{\fm})))\\
&\cong {\Hom}_{A}(N \otimes_{\bar{A}}\bar{A}, E(A/{\fm}))\\
&\cong {\Hom}_{A}(N, E(A/{\fm})) .\end{array}$$

\indent By Independence Theorem for the local cohomologies, we
have $H_{\bar{\widehat{{\fm}}}}^{s-i}(K_{\widehat{M}})\cong
H_{\widehat{{\fm}}}^{s-i}(K_{\widehat{M}})$.
 Put all these together, we obtain, from (4) and 1.2, that
 \begin{center}
 ${\Hom}_{A}(H^{i}(C_{A}(M)'), E(A/{\fm}))\cong
 H_{\widehat{{\fm}}}^{s-i}(K_{\widehat{M}}),
 i=0, 1, ..., $\end{center} as $A$ and $\widehat{A}$--modules.\\ \indent If
 $M$ admits a canonical module $K_{M}$, we then have
 $\widehat{(K_{M})} \cong K_{\widehat{M}}$, and by the Artinianness of
 $H_{{\fm}}^{s-i}(K_{M})$, we get the final claim. $\Box$
 \section {Applications}
\noindent First we show that over a local ring with
Cohen--Macaulay formal fibres, certain $f$--modules are also
generalized Cohen--Macaulay modules. Recall that $M$ is called
{\it generalized Cohen--Macaulay} (abbr. g.CM) if there exists
$r\geq1$ such that, for each system of parameters $x_{1}, ...,
x_{s}$ for $M$ and for all $i=1, ..., s$,  \begin{center}
${\fm}^{r}[((x_{1}, ..., x_{i-1})M : x_{i})/(x_{1}, ...,
x_{i-1})M]=0.$
\end{center} Note that, by [ScTC, (3.2) and (3.3)], $M$ is a g.CM
module if and only if $H_{{\fm}}^{i}(M)$ is of finite length for
all $i=0, 1, ..., s-1$ .\\ \indent An $A$--module $M$ is called an
$f$--module if for each system of parameters $x_{1}, ..., x_{s}$
for $M$ \begin{center} ${\Supp}_{A}[((x_{1}, ..., x_{i-1})M
:x_{i})/(x_{1}, ..., x_{i-1})M] \subseteq {\{{\fm}\}}$
\end{center} for all $i=1, ..., s$. It is clear that if $M$ is g.CM
module then it is an $f$--module.
\subsection{Theorem} (Compare [ScTC, (3.8)]).
{\it Assume that all formal fibres of A are Cohen--Macaulay. Let M
be an A--module such that }
${\Min}_{\widehat{A}}(\widehat{M})={\Assh}_{\widehat{A}}(\widehat{M})$.
{\it If M is an f--module with} ${\depth}_{A}(M)\geq2$,  {\it then
M is a
g.CM module}.\\

\noindent {\it Proof}. By a straightforward argument and using the
equivalent definition of $f$--module [T, Lemma 1.2 (ii)], it can
be shown that $M$ is $(S_{2})$ and that
${\Min}_{A}(M)={\Assh}_{A}(M)$.\\
\indent Now, by [T, Lemma 1.2 (iv)], the $M$--height filtration of
Spec($A$) is the same as the $M$--dimension filtration
$\mathcal{D}$ of Spec($A$), where $\mathcal{D}=(D_{i})_{i\geq0},
D_{i}=\{{\fp}\in{\Supp}_{A}(M): {\dim}(A/{\fp})\leq s-i\}$. Thus,
by [DT1, Lemma 3.1], there exists an isomorphism
\begin{center}
$C_{A}(M)=C(\mathcal{D}, M)\cong C(\mathcal{U}, M)$ \end{center}
 (over ${\Id}_{M}$), where $C(\mathcal{U}, M)$ is the complex of
 modules of generalized fractions on $M$ with respect to the chain
 of triangular subsets $\mathcal{U}=(U_{i})_{i\geq 1}$ on $A$,
 defined by \begin{center} $U_{i}=\{(x_{1}, \cdots, x_{i})\in A^{i} :
  $
 there exists $j$ with $0\leq j\leq i$ such that $x_{1}, \cdots,
 x_{j}$\\
 is an s.s.o.p. for $M$ and $x_{j+1}=\cdots=x_{i}=1\}$
 \end{center} (See [DT1] for details).  By [SZ, Corollary 2.3 and
 Theorem 2.4], \begin{center} $H^{i-1}(C_{A}(M))\cong
 H_{{\fm}}^{i}(M), i=1, \cdots, s-1$. \end{center} Therefore, by
 Theorem 2.1, $H_{{\fm}}^{i}(M)$ is of finite length for all
 $i=0, 1, \cdots, s-1.  \Box$
  \subsection{Corollary} {\it Assume
 that all formal fibres of A are Cohen-Macaulay. If A is an f-ring
 with} ${\depth}(A)\geq 2$, {\it then A is a g.CM ring}.\\

 \noindent {\it Proof}. As we have seen in the proof of 3.1, $A$
 is $(S_{2})$. By [M, Theorem 23.9], $\widehat{A}$ satisfies $(S_{2})$.
 Thus ${\Min}(\widehat{A})={\Assh}(\widehat{A})$ (see [AG, Lemma 1.1]). Now
 the result follows from Theorem 3.1. $\Box$ \\\
 \indent Our next
 application studies the injectivity of the terms of the Cousin
 complex $C_{A}(M)$.\\\\
 \indent In [S2],
 a finitely generated $A$-module $M$ is defined to be a Gorenstein $A$-module
 whenever its Cousin complex provides a minimal injective
 resolution. It is also proved that if $A$ admits a canonical
 module $\Omega$, then any Gorenstein $A$-module is isomorphic to
 the direct sum of a finite number of copies of $\Omega$  [S3, Theorem
 2.1].\\
  \indent It is known that if $A$ does not have a canonical
 module and has a Gorenstein module, then it has a unique
 indecomposable Gorenstein module $G$ and every Gorenstein
 $A$-module is isomorphic to a direct sum of a finite number of
 copies $G$ (see[FFGR and S2]). Here we extend this result and show
 that for any finitely generated module $M$, over a complete $(S_{2})$
 local ring $A$ which satisfies $(S_{2})$, if $0:_{A}M=0$ and
 $C_{A}(M)'$ is an injective complex, then $M$ is isomorphic to a
 direct sum of copies of a uniquely determined indecomposable
 one.\\
 \subsection{Theorem} {\it Let A satisfy} $(S_{2})$ {\it and
 suppose that it possesses a dualizing complex. Assume that M
 satisfies} $(S_{2})$ {\it and} $0:_{A}M=0$. {\it The following
 statements are equivalent}:\\ \indent (i) $C_{A}(M)'$ {\it is an injective
 complex};\\ \indent (ii) $M$ {\it is isomorphic to a direct sum
 of a finite number of copies of the canonical module $K$ of the ring
 } $A$.\\

 \noindent {\it Proof}. (i)$\Rightarrow$ (ii). We do not need
 $A$ to satisfy $(S_{2})$ in this part. The proof is a
 straightforward adaptation of the argument in [S3, Theorem
 2.1(v)]. Let $K$ denote the canonical module of $A$.\\ \indent
 Let \begin{center} $I^{\bullet} : 0\longrightarrow I^{0}\longrightarrow
  I^{1}\longrightarrow \cdots\longrightarrow I^{d}\longrightarrow
  0$ \end{center}
  be the dualizing complex for $A$ so that $K=H^{0}(I^{\bullet})$. By
  the proof of [DT2, Lemma 3.1], $C_{A}(M)\cong
  {\Hom}_{A}(K_{M}, I^{\bullet})^\ast$, where
  $K_{M}={\Hom}_{A}(M, K)$. Hence all cohomology modules
  of $C_{A}(M)$ are finitely generated (see  [S4, Lemma 3.4(ii)]).
  By [S5, Theorem], ${\Hom}_{A}(K_{M}, I^{d}) \cong
  H_{{\fm}}^{d}(M)$. As  $H_{{\fm}}^{d}(M)$ is an Artinian
  injective $A$-module, we may write $H_{\fm}^{d}(M) \cong
  \oplus_{i=1}^{n}E(A/{\fm})$, say. Using the Matlis functor
  ${\Hom}_{A}(- , E(A/{\fm}))$ and that $I^{d}=E(A/{\fm})$, we
  obtain $K_{M}\otimes_{A}\widehat{A}\cong
  (\oplus_{i=1}^{n}A)\otimes_{A}\widehat{A}$. This implies, by[HK,
  Lemma 5.8], that $K_{M}\cong \oplus_{i=1}^{n}A$. Hence we have
  $H_{{\fm}}^{d}(K_{M}) \cong \oplus_{i=1}^{n}H_{{\fm}}^{d}(A)$.\\
  \indent On the other hand, by Grothendieck's Local Duality
  Theorem [B-ZS, Corollary 2.5] and the fact that $M$ satisfies
  $(S_{2})$ so $C_{A}(M)$ is exact at point $-1, 0$ (see [SSc,
  Example 4.4]), we obtain\ \begin{center} $H_{{\fm}}^{d}(K_{M}) \cong
  {\Hom}_{A}(H^{0}(C_{A}(M)'), E(A/{\fm})) \cong {\Hom}_{A}(M,
  E(A/{\fm}))$. \end{center}
  By applying the Matlis functor again, we get
  $M\otimes_{A}\widehat{A} \cong {\Hom}_{A}(H_{{\fm}}^{d}(K_{M}),
  E(A/{\fm})) \cong {\Hom}_{A}(\oplus_{i=1}^{n}H_{{\fm}}^{d}(A),
  E(A/{\fm})) \cong (\oplus_{i=1}^{n}K)\otimes_{A}\widehat{A}$. Now,
  by [HK, Lemma 5.8], $M\cong \oplus_{i=1}^{n}K$.\\
  \indent (ii)$\Rightarrow$(i). We have
  ${\Supp}_{A}(M)={\Supp}_{A}(K)= {\Spec}(A)$ (see [A, (1.8)] and
  [AG, Lemma 1.1]). It is routine to check that $C_{A}(M) \cong
  \oplus_{i=1}^{n}C_{A}(K)$. As ${\Min}A={\Assh}A$ and the
  dimension filtration and the height filtration of ${\Spec}(A)$
  are the same (see [A, (1.9)]), the claim follows by [DT1,
  Corollary 3.4].$\Box$
  \subsection{Corollary} {\it Assume that} $\widehat{A}$ {\it
  satisfies}
  $(S_{2})$. {\it Then the following statements are equivalent} :
  \\
  \indent (i) $C_{\widehat{A}}(\widehat{A})'$ {\it is an injective complex of
  $\widehat{A}$-modules;} \\  \indent  (ii) {\it A is the canonical
  module of A.} \\  \indent {\it Moreover, if A satisfies one of the
  above equivalent conditions, then A is Gorenstein if and only if
  $\widehat{A}$ satisfies $(S_{n})$ for some} $n\geq (1/2){\dim}A+1$.\\

   \noindent {\it Proof}. (i)$\Rightarrow$ (ii). Set $\Omega$ for
   the canonical module of $\widehat{A}$. By 3.3, $\widehat{A} \cong
   \Omega^{n}$ for some $n$. Thus $H_{\widehat{{\fm}}}^{d}(\widehat{A})
   \cong\oplus_{i=1}^{n}H_{\widehat{{\fm}}}^{d}(\Omega) \cong
   \oplus_{i=1}^{n}E(\widehat{A}/\widehat{{\fm}})$ and, by applying
   ${\Hom}_{\widehat{A}}(-, E(\widehat{A}/\widehat{{\fm}}))$, we get $\Omega
   \cong \widehat{A}^{n}$. Thus $\widehat{A}^{n^{2}}=\widehat{A}$, which
   implies $n=1$ and so $A$ is the canonical module of $A$. \\
   \indent (ii)$\Rightarrow$(i). As $\widehat{A}$ is the canonical
   module of $\widehat{A}$ and $\widehat{A}$ satisfies $(S_{2})$,
   $C_{\widehat{A}}(\widehat{A})'$ is the dualizing complex of $\widehat{A}$
   [DT1, Corollary 3.4]. \\ \indent For the last part, we may
   assume that $A$ is complete. By [SSc, Example 4.4], $C_{A}(A)$
   is exact at points $-1, 0, 1, ..., n-2$,  from which it follows, by
   Theorem 2.4, that $H_{{\fm}}^{{\dim}A-i}(A)=0$ for
   $0<i\leq n-2$. On the other hand, as $A$ satisfies $(S_{n})$,
   $H_{{\fm}}^{i}(A)=0$ for all $i< $min$\{d, n\}$. As
   ${\dim}A-(n-2)\leq n$, it follows that $H_{{\fm}}^{i}(A)=0$ for
   all $i<d$, which imply the exactness of $C_{A}(A)$. The other
   side is trivial.$\Box$ \\
   \section{Canonical modules of non--local rings}
   Recall that, for a Noetherian (not necessarily local) ring $R$, the canonical module
   of R (if it exists) is a finite $R$-- module $K$ such that
   $K_{{\fm}}$, the localization of $K$ at any maximal ideal
   ${\fm}$ of $R$, is the canonical module of $R_{{\fm}}$. In
   order to generalize our results to the non--local case one might
   ask whether a canonical module exists even when $R$ possesses
   a dualizing complex. We will show that, if $R$ satisfies
   $(S_{2})$ and all formal fibres of $R_{{\fm}}$, for any maximal ideal
   ${\fm}$ of $R$, are
   Cohen-Macaulay, then existence of a canonical module for $R$ is
   equivalent to the statement that $R$ possesses a dualizing complex.\\ \indent
   Throughout, $R$ is a Noetherian ring of finite dimension
   which is not necessarily  local.\\
   \indent Assume that $R$ possesses a dualizing complex
   $I^{\bullet}$ and $t({\fp}; I^{\bullet}), {\fp}\in
   {\Spec}R$, denotes the unique integer $i$
   for which ${\fp}$ occurs in $I^{i}$ (see [H, page 23]).\\
   \subsection{Proposition} {\it Assume that R satisfies
   $(S_{2})$ and that it possesses a dualizing complex
   $I^{\bullet}$. If ${\fp}, {\fq} \in {\Min}(R)$ such that
   ${\fp} \subseteq {\fr}$ and ${\fq}\subseteq {\fr}$ for some
   ${\fr}\in {\Spec}(R)$, then} $t({\fp}; I^{\bullet})= t({\fq};
   I^{\bullet})$.\\

   \noindent {\it Proof}. We may assume that $R$ is
   a local ring and that its maximal ideal is ${\fr}$. As $R$
   satisfies $(S_{2})$ and possesses a dualizing complex, then
   ${\Min}(R) ={\Assh}(R)$ [A; 1.1]. Therefore $t({\fp}; I^{\bullet})= t({\fq};
   I^{\bullet}).  \Box$
   \subsection{Notation} Assume that $R$ satisfies $(S_{2})$ and
   that \\ \begin{center} $I^{\bullet}: 0\longrightarrow
   I^{0}\overset{\delta^{0}}
   \longrightarrow I^{1}\overset{\delta^{1}}\longrightarrow\cdots
   \overset{\delta^{l-1}}
   \longrightarrow I^{l}\longrightarrow 0$ \end{center} is a
   dualizing complex for $R$. It follows that ${\Ass}_{R}(I^{0})\subseteq
   {\Min}R$.  Assume that
   ${\Ass}_{R}(I^{0})\neq{\Min}(R)$.  Let $r$ be the greatest
   integer such that $X :=
   {\Min}(R)\bigcap{\Ass}_{R}(I^{r})\neq \emptyset$. Set, for each
   $i\geq 0$,  \begin{center}
   $X_{i} = \{{\fp}\in {\Ass}_{R}(I^{i}): {\fp}$ contains some element of
   $X\}$; \end{center} \begin{center} $X_{i}'={\Ass}_{R}(I^{i}) \setminus X_{i}$;
   \end{center} \begin{center} $I_{1}^{i} =
   \oplus_{{\fp}\in X_{i}}E(A/{\fp}),  I_{2}^{i} =
   \oplus_{{\fp}\in X_{i}'}E(A/{\fp})$, \end{center} \noindent so
   that we have $I^{i} = I_{1}^{i}\oplus I_{2}^{i}$.\\
   \subsection{Proposition} {\it With the notations as in 4.2}, \
   \begin{center}
   ${\Hom}_{R}(I_{1}^{i}, I_{2}^{i+1}) = 0 = {\Hom}_{R}(I_{2}^{i},
   I_{1}^{i+1})$.\end{center}

   \noindent {\it Proof}. If ${\Hom}_{R}(I_{1}^{i},
    I_{2}^{i+1}) \neq 0$, then ${\Hom}_{R}(I_{1}^{i},
    E(R/{\fp}))\neq 0$ for some ${\fp}\in X_{i+1}'$. Assume that
    $f : I_{1}^{i}\longrightarrow E(R/{\fp})$ is an
    $R$--homomorphism and that $x\in I_{1}^{i}$. Let
    $x=x_{1}+\cdots+x_{s}$, where $x_{j}\in E(R/{\fp}_{j}),
    1\leq j\leq s$. By definition of $X_{i+1}'$, we have
    ${\fp}_{1}\bigcap\cdots \bigcap {\fp}_{s}\nsubseteqq {\fp}$.
    Take $t\in {\fp}_{1}\bigcap\cdots \bigcap {\fp}_{s} \setminus {\fp}$.
    Hence $t^{m}x=0$ for some positive integer $m$. On the other
    hand the map $E(R/{\fp})\overset{t^{m}}\longrightarrow
    E(R/{\fp})$ is an isomorphism. Thus $t^{m}f(x)=f(t^{m}x)=0$
    implies that $f(x)=0$.\\
    \indent To show that ${\Hom}_{R}(I_{2}^{i}, I_{1}^{i+1})= 0$,
    we may assume, on the contrary, that
    ${\Hom}_{R}(I_{2}^{i}, E(R/{\fp})) \neq 0$ for some ${\fp}\in
    X_{i+1}$. So we may assume that
    ${\fp}'\subseteq {\fp}$, for some ${\fp}'\in X_{i}'$. By
    localizing at ${\fp}$, we get ${\Min}(R_{{\fp}})=
    {\Assh}(R_{{\fp}})$, because $R_{{\fp}}$ satisfies $(S_{2})$
    and $R_{{\fp}}$  possesses a dualizing complex. As
    ${\fp}'R_{{\fp}}$  contains a minimal element
    ${\fq}R_{{\fp}}\in{\Assh}(R_{{\fp}})$, we get, by 4.1,
    $t({\fq}; I^{\bullet})= r$ and that ${\fq}\in X$. This contradicts with the
    definition of $X_{i}'.\Box$
   \subsection{Theorem} {\it Assume that R satisfies $(S_{2})$ and
   that it possesses a dualizing complex. Then R possesses a
   dualizing complex \ \begin{center} $J^{\bullet} :
   0\longrightarrow J^{0}\longrightarrow J^{1}\longrightarrow
   \cdots \longrightarrow J^{d}\longrightarrow 0,    d={\dim}R,$
   \end{center} such that ${\Ass}_{R}(J^{0})= {\Min}(R)$. In
   particular $R$ admits a canonical module.}\\

   \noindent {\it Proof}. The proof is influenced by  [H, Lemma 3.1]. Suppose that\
     \begin{center} $I^{\bullet}:0\longrightarrow
     I^{0}\overset{\delta^{0}}\longrightarrow I^{1}
     \overset{\delta^{1}}\longrightarrow\cdots \overset{\delta^{l-1}}
     \longrightarrow I^{l}\longrightarrow 0$ \end{center}
is the dualizing complex for $R$. Assume further that
${\Ass}(I^{0})\neq{\Min}(R)$ and that $r$ is the greatest integer
with $X:={\Min}(R)\bigcap{\Ass}_{R}(I^{r})\neq\emptyset$. We set
$X_{i}, X'_{i}, I_{1}^{i}$ and $I_{2}^{i}$ as in 4.2. Note that
$I_{1}^{i}=0$ and $I_{2}^{i}=I^{i}$ for $0\leq i<r$
(see [S6, Lemma 3.3]).\\
\indent We construct a dualizing complex \
\begin{center} $J^{\bullet}: 0\longrightarrow
J^{0} \overset{\eta_{0}}\longrightarrow J^{1}
\overset{\eta_{1}}\longrightarrow \cdots$ \end{center} as follows.
Set $J^{i}=I_{1}^{i+r}\oplus I_{2}^{i}$ for all $i\geq 0$, and
define $\eta^{i}:J^{i}\longrightarrow J^{i+1}$ by
$\eta^{i}(x+y)=\delta_{1}^{i+r}(x)+\delta_{2}^{i}(y)$ for $x\in
I_{1}^{i+r}, y\in I_{2}^{i}$, where
$\delta_{1}^{j}:=\delta^{j}\mid_{I_{1}^{j}}$ and $
\delta_{2}^{j}:=\delta^{j}\mid_{I_{2}^{j}}, j\geq 0$. It follows
from Proposition 4.3 that $J^{\bullet}$ is a complex. To show that
$H^{i}(J^{\bullet})$ is a finitely generated $R$--module for all
$i\geq 0$, we note that, by a straightforward argument, there are
two natural isomorphisms \\ \begin{center}
$H^{i}(I^{\bullet})\cong({\Ker}\delta_{1}^{i}/{\im}\delta_{1}^{i-1})
\oplus({\Ker}\delta_{2}^{i}/{\im}\delta_{2}^{i-1})$, \end{center}
\begin{center}
$H^{i}(J^{\bullet})\cong({\Ker}\delta_{1}^{i+r}/{\im}\delta_{1}^{i+r-1})
\oplus({\Ker}\delta_{2}^{i}/{\im}\delta_{2}^{i-1}), i\geq 0.$
\end{center}
Therefore $J^{\bullet}$ is a dualizing complex for $A$.\ \indent
Now we have $J^{0}=I_{1}^{r}\oplus I^{0}$ and thus
${\Ass}_{R}(I^{0})\subsetneqq{\Ass}_{R}(J^{0})$. So after a finite
number of steps we are finished.\\ \indent Finally, let
$J^{\bullet}$ be a dualizing complex with
${\Min}(R)={\Ass}_{R}(J^{0})$. For each ${\fm}\in{\Max}(R)$, the
complex\
\begin{center}
$0\longrightarrow(J^{0})_{{\fm}}\longrightarrow(J^{1})_{{\fm}}
\longrightarrow\cdots\longrightarrow(J^{t({\fm};J^{\bullet}}))_{{\fm}}
\longrightarrow 0$ \end{center} is the dualizing complex for
$R_{{\fm}}$, so that, by Grothendieck's Local Duality Theorem
[B-ZS, Corollary 2.5], $(H^{0}(J^{\bullet}))_{{\fm}}$ is the
canonical module of
$R_{{\fm}}$. Thus $H^{0}(J^{\bullet})$ is a canonical module of $R$. $\Box$\\
\indent As an application, we can give a partial generalization of
[BH, Proposition 3.3.18].\\
\subsection{Theorem} {\it Assume that R
satisfies $(S_{2})$ and possesses a dualizing complex \
\begin{center} $I^{\bullet}:0\longrightarrow I^{0}\longrightarrow
I^{1}\longrightarrow \cdots\longrightarrow I^{d}\longrightarrow
0$, \end{center} with $d={\dim}R,
I^{i}=\underset{{\heit}{\fp}=i}\oplus E(R/{\fp}), i=0, 1, \cdots$
 and that $K_{R}$ is a canonical module of $R$.\\
\indent (a)  The following conditions
are equivalent:\\ \indent  \indent  (i) $K_{R}$ has a rank;\\
\indent
 \indent (ii) ${\rank}K_{R}=1$;\\ \indent  \indent (iii) $R$ is
generically Gorenstein ( that is $R_{{\fp}}$ is a Gorenstein ring
for all mini-\\ \indent \indent \indent \indent mal prime ideals ${\fp}$ of $R$).\\
\indent (b) If $K_{R}$ satisfies $(S_{3})$ and the equivalent
conditions of (a) hold, then $K_{R}$ can be identified with an
ideal of height 1 or equals $R$. In the first case $R/K_{R}$, is
an $(S_{2})$ ring with the canonical module  $R/K_{R}$, the ring
itself.}\\

\noindent {\it Proof}. The proof is parallel to that of [BH,
Proposition
3.3.18] and we present it for the convenience of the reader.\\
\indent (a). (i)$\Rightarrow$(ii)$\Rightarrow$(iii). Set $Q$ to be
the ring of total fractions of $R$ and let $K_{R}\otimes_{R}Q$ be
a free $Q$--module of rank $r$, say. Let ${\fp}\in{\Min}(R)$. As
${\Min}(R)= {\Ass}(R)= {\Ass}_{R}(K_{R})$ (see [DT1, 1.3]), we
know that $K_{R_{{\fp}}}\cong (K_{R})_{{\fp}}$ and, by [BH,
Proposition 1.4.3], the $R_{{\fp}}$--module $(K_{R})_{{\fp}}$ is
free of rank $r$. As the dualizing complex of $R_{{\fp}}$ is
$(I^{\bullet})_{{\fp}}: 0\longrightarrow
I_{{\fp}}^{0}\longrightarrow 0$, we get $(R_{{\fp}})^{r}\cong
(K_{R})_{{\fp}}\cong E(R/{\fp})$ from which
it follows that $r=1$, and thus $R_{{\fp}}$ is Gorenstein.\\
\indent \indent (iii)$\Rightarrow$(i). Note that
${\Min}(R)={\Ass}(R)$, and thus [BH, Proposition 1.4.3] implies
that $K_{R}$ has rank 1.\\ \indent As ${\Ass}R= {\Ass}(K_{R}), K$
is torsion free. Thus [BH, 1.4.18] implies that $K_{R}$ is
isomorphic to a sub--module of a free $R$--module of rank 1, and
it may be identified with an ideal of $R$ which we again denote by
$K_{R}$.\\ \indent If ${\dim}R=0$, we get $K_{R}\cong R$, so we
may assume ${\dim}R>0$, and also $K_{R}$ is a proper ideal of $R$.
By [BH, Proposition 1.4.3], $K_{R}$ has a free sub--module
${\fa}$, which is also an ideal of $R$ of rank 1. Assuming ${\fa}=
xR$ with $x$ is a base for ${\fa},  x$ is $R$--regular and
$K_{R}$--regular.\ \indent Let ${\fp}$ be a prime ideal containing
$K_{R}$. Applying the functor ${\Hom}_{R_{{\fp}}}(- ,
(I^{\bullet})_{{\fp}})$ on the exact sequence $0\longrightarrow
K_{R}R_{{\fp}}\longrightarrow R_{{\fp}}\longrightarrow
R_{{\fp}}/K_{R}R_{{\fp}}\longrightarrow 0$, we get the exact
sequence
$$\begin{array}{lll} 0 &\longrightarrow
H^{0}({\Hom}_{R_{{\fp}}}(R_{{\fp}}/K_{R}R_{{\fp}},
(I^{\bullet})_{{\fp}}))\longrightarrow
H^{0}((I^{\bullet})_{{\fp}})\longrightarrow
H^{0}({\Hom}_{R_{{\fp}}}(K_{R}R_{{\fp}}, (I^{\bullet})_{{\fp}}))\\
&\longrightarrow
H^{1}({\Hom}_{R_{{\fp}}}(R_{{\fp}}/K_{R}R_{{\fp}},
(I^{\bullet})_{{\fp}}))\longrightarrow
H^{1}((I^{\bullet})_{{\fp}})\longrightarrow \cdots.\end{array}$$

Note that $H^{1}(I^{\bullet})_{{\fp}}=0$  as $K_{R_{{\fp}}}$
satisfies $(S_{3})$ (see [DT2, Proposition 2.5]).\\ \indent On the
other hand, we have
$H^{0}({\Hom}_{R_{{\fp}}}(R_{{\fp}}/K_{R}R_{{\fp}},
(I^{\bullet})_{{\fp}})) \cong 0:_{K_{R}R_{{\fp}}}
K_{R}R_{{\fp}}\subseteq 0:_{R_{{\fp}}} K_{R}R_{{\fp}}= 0$ because
$K_{R}R_{{\fp}}$ is the canonical module of $R_{{\fp}}$ and
$R_{{\fp}}$ satisfies $(S_{2})$. Therefore we get the exact
sequence \ \begin{center} $o\longrightarrow
K_{R}R_{{\fp}}\longrightarrow R_{{\fp}}\longrightarrow
H^{1}({\Hom}_{R_{{\fp}}}(R_{{\fp}}/K_{R}R_{{\fp}},
(I^{\bullet})_{{\fp}}))\longrightarrow 0$ \end{center} which
implies that $H^{1}({\Hom}_{R_{{\fp}}}(R_{{\fp}}/K_{R}R_{{\fp}},
(I^{\bullet})_{{\fp}}))\cong (R/K_{R})_{{\fp}}$. It follows that
${\dim}(R_{{\fp}}/K_{R}R_{{\fp}})= {\heit}{\fp}-1$. Using the
Grothendieck local duality shows that $R_{{\fp}}/K_{R}R_{{\fp}}$
is the canonical module of $ R_{{\fp}}/K_{R}R_{{\fp}}$. As $
K_{R}R_{{\fp}}$ contains an $R_{{\fp}}$--regular element, we have
${\heit}_{R_{{\fp}}}( K_{R}R_{{\fp}})\geq 1$. Since
${\dim}(R_{{\fp}}/K_{R}R_{{\fp}})= {\heit}{\fp}-1$, we get
${\heit}(K_{R})=1$.\\
\indent For the final part, we may assume that ${\dim}R> 3$. As
$R$ is $(S_{2})$ and $K_{R}$ is $(S_{3})$, from the exact sequence
$H_{{\fp}R_{{\fp}}}^{i}(R_{{\fp}})\longrightarrow
H_{{\fp}R_{{\fp}}}^{i}(R_{{\fp}}/K_{R}R_{{\fp}})\longrightarrow
H_{{\fp}R_{{\fp}}}^{i+1}(K_{R}R_{{\fp}})$, we get $
H_{{\fp}R_{{\fp}}}^{i}(R_{{\fp}}/K_{R}R_{{\fp}})= 0$ for $i=0, 1$.
This shows that $R_{{\fp}}/K_{R}R_{{\fp}}$ satisfies
$(S_{2}).\Box$\\\
 \indent We can also generalize[DT1, Corollary 3.4].\\
 \subsection{Theorem} {\it Assume that $R$ satisfies $(S_{2})$, and
 that ${\dim}R<\infty$. The following statements are equivalent.\\
 \indent (i) R possesses a dualizing complex;\\
 \indent (ii) R admits a canonical module K, and $C(\mathcal{H},
 K)'$, the induced complex of the Cousin complex of K with respect
 to the height filtration $\mathcal{H}=(H_{i})_{i\geq0}$ with
 $H_{i}=\{{\fp}\in {\Spec}(R): {\heit}({\fp})\geq i\}$, is a
 dualizing complex for R;\\
 \indent (iii)R admits a canonical module K and $H^{i}(C(\mathcal{H},
 K)')$ is finitely a generated R--module for all $i\geq1$.}\\

\noindent {\it Proof}. (i)$\Rightarrow$(ii). By 4.4, there exists
a
 dualizing complex\ \begin{center} $I^{\bullet}: 0\longrightarrow
 I^{0} \overset{\delta^{0}}\longrightarrow I^{1}
 \overset{\delta^{1}}\longrightarrow \cdots$ \end{center} for $R$ such that
 ${\Ass}_{R}(I^{0})= {\Min}(R)$. Set $K={\Ker}\delta^{0}$. As seen
 in 4.4, $K$ is a canonical module of $R$ and
 ${\Ass}_{R}(K)={\Min}(R)$. Now, by [DT1, Theorem 2.4(iv)],
 $C(\mathcal{H}, K)'$ is a dualizing complex for $R$.\\ \indent
 (ii)$\Rightarrow$(iii) is clear.\\ \indent (iii)$\Rightarrow$(i).
 For each ${\fm}\in {\Max}(R)$, we have, by [S1, Theorem 3.5],
  $C(\mathcal{H}, K)_{{\fm}}\cong C(\mathcal{H_{{\fm}}},
  K_{{\fm}})$, where $\mathcal{H}_{{\fm}}$ is the height filtration
  of $R_{{\fm}}$. Therefore, by [DT1, Corollary 3.4], $(C(\mathcal{H}
  , K)')_{{\fm}}$ is a dualizing complex for $R_{{\fm}}$. Since $R_{{\fp}}$
  satisfies $(S_{2})$ for all ${\fp}\in {\Spec}(R)$, by the same
  argument as in the proof of [DT1, Corollary 3.4], each term of
  $C(\mathcal{H}, K)'$ is an injective module. Thus, by [S4,
  Theorem 4.2], $C(\mathcal{H}, K)'$ is a dualizing complex for
  $R$. $\Box$\\
  \section{Indecomposable injective modules structure} In this
  section, by using of a particular dualizing complex for an
  $(S_{2})$ Noetherian ring $R$ of finite dimension, we give an explicit description for the
  structure of all indecomposable injective modules. In [DT1,
  Corollary 3.3], it is shown that for each ${\fp}\in{\Spec}(R)$,
  there exists a finitely generated $R$--module $T$, depending on
  ${\fp}$, such that $E(R/{\fp})$ is a module of generalized
  fractions of $T$. Here we will show that $T$ can be replaced  by a canonical
  module of $R$ and that it does not depend on ${\fp}$.\\\
  \indent Our approach involves the concept of a chain of
  triangular subsets on $R$ explained in [O, page 420]. Such a
  chain $\mathcal{U}=(U_{i})_{i\geq 1}$ determines a complex
  $C(\mathcal{U}, M)$ of modules of generalized fractions on an
  $R$--module $M$, that is\\ \begin{center} $C(\mathcal{U}, M):
  0\longrightarrow M \overset{e^{0}}\longrightarrow U_{1}^{-1}M
  \overset{e^{1}}\longrightarrow \cdots\overset{e^{i-1}}\longrightarrow U_{i}^{-i}M
  \overset{e^{i}}\longrightarrow
  U_{i+1}^{-i-1}M\overset{e^{i+1}}\longrightarrow\cdots$ \end{center} in which
  $e^{0}(m)=\frac{m}{(1)}$ for all $m\in M$ and
  $e^{i}(\frac{m}{(u_{1}, \cdots, u_{i})})=
  \frac{m}{(u_{1}, \cdots, u_{i}, 1)}$ for all $i\geq1$, $m\in M$,
  and $(u_{1}, \cdots, u_{i})\in U_{i}$. Note that in the complex
  $C(\mathcal{U},
  M)$, $U_{i+1}^{-i-1}M$ is regarded as the $i$--th term, so that
  $H^{i}(C(\mathcal{U}, M))= {\Ker}e^{i+1}/{\im}e^{i}$, $i\geq 0$,

  and $H^{-1}(C(\mathcal{U}, M))= {\Ker}e^{0}$.\\
  \indent Assume that $R$ satisfies $(S_{2})$ and possesses a
  dualizing complex, so that $R$ possesses a dualizing complex \
  \begin{center} $I^{\bullet}:0\longrightarrow
  I^{0}\overset{\delta^{0}}\longrightarrow I^{1} \overset{\delta^{1}}
  \longrightarrow \cdots\overset{\delta^{d-1}}\longrightarrow I^{d}\longrightarrow 0,
  d={\dim}R$, \end{center}
  such that ${\Ass}_{R}(I^{0})={\Min}(R)$. Set
  $K={\Ker}\delta^{0}$, and consider the induced extended complex

\begin{center} $I^{\ast}:0\longrightarrow K\hookrightarrow
  I^{0}\overset{\delta^{0}}\longrightarrow I^{1} \overset{\delta^{1}}
  \longrightarrow \cdots\longrightarrow I^{d}\longrightarrow 0, $ . \end{center}
   For each ${\fp}\in{\Ass}_{R}(I^{0})$, the complex
$0\longrightarrow (I^{0})_{{\fp}}\longrightarrow 0$ is the
dualizing complex for $R_{{\fp}}$, so that $K_{{\fp}}\cong
E(R/{\fp})$. Hence ${\Ass}_{R}(K)= {\Min}(R)$. Thus, by [DT1,
Proposition 3.2], there is a unique isomorphism of complexes (over
${\Id}_{K}$) from $I^{\ast}$ to $C(\mathcal{V}, K)$, the complex
of modules of generalized fractions on $K$ with respect to the
chain of triangular subsets $\mathcal{V}=(V_{i})_{i\geq 1}$ on
$R$, defined by\
 \begin{center} $V_{i}= \{(v_{1}, \cdots, v_{i})\in R^{i}:
 {\heit}_{R}((v_{1}, \cdots, v_{j}))\geq j$ for all $j$ with $1\leq
 j\leq i\}$. \end{center}
 \indent Now, we restate [DT1, Corollary 3.3] in a more
 appropriate form.
 \subsection{Corollary} {\it Assume that $R$ satisfies $(S_{2})$
 and that it possesses a dualizing complex, so that $R$ admits a
 canonical module $K$, say. Then, for each ${\fp}\in {\Spec}(R)$,
 \
  \begin{center} $E(R/{\fp})\cong (V_{{\heit}{\fp}}\times
  (R\setminus{\fp}))^{-{\heit}{\fp}-1}K$ \end{center} where $V_{r}$
  is the triangular subset of $R^{r}$ defined in the paragraph
  just before the corollary.}\\

   \noindent {\bf Acknowledgment}. I thank M. Tousi for his comment on
   2.3. I also thank the referee for the invaluable comments on the
   manuscript.


\begin{thebibliography}{9999}



\bibitem[A]{1} {\sc Y. Aoyama, } {\sl Some basic results on canonical modules},
J. Math. Kyoto Univ. {\bf 23} (1983), 85--94.

\bibitem[AG]{2} {\sc Y. Aoyama and S. Goto, } {\sl On the endomorphism ring
 of the canonical module}, J. Math. Kyoto Univ. {\bf 25} (1985),
 21--30.

\bibitem[B-ZS]{4} {\sc M. H. Bijan-Zadeh and R. Y. Sharp, }
{\sl On Grothendieck's local duality theorem},  Math. Proc.
Cambridge Philos. Soc. {\bf 85} (1979),  431--437.

\bibitem[BH]{3} {\sc W. Bruns and J. Herzog, } {\sl Cohen--Macaulay
Rings, } Cambridge University Press,  1996.

\bibitem[DT1]{5} {\sc M. T. Dibaei and M. Tousi, }
{\sl The structure of dualizing complex for a ring which is
$(S_{2})$}, J. Math. Kyoto Univ. {\bf 38} (1998),  503--516.

\bibitem[DT2]{6} {\sc M. T. Dibaei and M. Tousi, } {\sl A generalization of dualizing
complex structure and its applications}, J. Pure and Applied
Algebra,  {\bf 155} (2001), 17--28.



\bibitem[FFGR]{7} {\sc R. Fossum, H.-B. Foxby, P. Griffith, and I. Reiten, } {\sl
Minimal injective resolutions with applications to dualizing
modules and Gorenstein modules, }  Inst. Hautes Etudes Sci. Publ.
Math.,  {\bf 45} (1976), 193--215.

\bibitem[H]{8} {\sc J. E. Hall,} {\sl Fundamental dualizing complexes for
 commutative Noetherian rings}, Quart. J. Math. Oxford
 {\bf 165} (1979), 21--32.

\bibitem[HK]{9} {\sc J. Herzog and E. Kunz,}  {\sl Der Kaninische
Modul eines Cohen--Macaulay Rings},
Lecture Notes Math. 238, Springer--Verlag, 1971.
\bibitem[K] {} {\sc T. Kawasaki}, {\sl Finiteness of Cousin
homologies}, preprint.

\bibitem[M]{10} {\sc H. Matsumura,} {\sl Commutative ring theory},
Cambridge University Press, 1992.

\bibitem[O]{11} {\sc L. O'Carrol, } {\sl On the generalized
fractions of Sharp and Zakeri}, J. Lodon Math. Soc. {\bf 28}
(1983) 417--427.

\bibitem[P]{12} {\sc H. Petzl,} {\sl Cousin complexes and flat ring
extentions,} Comm. Algebra, {\bf 25} (1997), 311--339.

\bibitem[ScTC]{13} {\sc P. Schenzel, N. V. Trung and N. T. Cuong,}
{\sl Verallgemeinerte Cohen-Macaulay-Moduln,} Math. Nachr. {\bf
85} (1978), 57--73.

\bibitem[S1]{14} {\sc R. Y. Sharp,} {\sl The Cousin complex for a module
 over a commutative Noetherian ring},
 Math. Z., {\bf 112} (1969), 340--356.

\bibitem[S2]{15} {\sc R. Y. Sharp,} {\sl Gorenstein modules},
 Math. Z. {\bf 115} (1970), 117--139.

\bibitem[S3]{16} {\sc R. Y. Sharp,} {\sl Finitely generated modules of
finite injective dimension over certain Cohen--Macualay rings,}
London math. Soc. {\bf 25} (1972), 303--328.

\bibitem[S4]{17} {\sc R. Y. Sharp,} {\sl Dualizing complexes for
commutative Noetherian ring}, Math. Proc. Cambridge Philos. Soc.
{\bf 78} (1975),  369--386.

\bibitem[S5]{18} {\sc R. Y. Sharp,} {\sl Local cohomology and the Cousin
complex for a commutative Noetherian ring,}  Math. Z., {\bf 153}
(1977), 19--22.

\bibitem[S6]{19} {\sc R. Y. Sharp,} {\sl A commutative Noetherian
ring which possesses a dualizing complex is acceptable,} Math.
Proc. Camb. Philos. Soc. {\bf 82} (1977), 197--213.

\bibitem[SSc] {20} {\sc R. Y. Sharp and P. Schenzel,} {\sl Cousin
complex and generalized Hughes complexes,} Proc. London Math.
Soc., {\bf 68} (1994), 499--517.

\bibitem[SZ] {21} {\sc R. Y. Sharp and H. Zakeri,} {\sl Generalized
fractions, Buchsbaum modules and generalized Cohen--Macaulay
modules,} Math. Proc. Cambridge Philos. Soc., {\bf 98} (1985),
429--436.

\bibitem[T] {22} {\sc N. V. Trung,} {\sl Toward a theory of
generalized Cohen--Macauly modules,} Nagoya Math. J. {\bf 102}
(1986), 1--49.

\end{thebibliography}
\end{document}